\documentclass[11pt, oneside]{article}   	
\usepackage{geometry}                		
\usepackage{graphicx}				
\usepackage{amssymb}

\usepackage{theorem}
\usepackage[colorlinks=true,citecolor=blue]{hyperref}
\theoremheaderfont{\normalfont\bfseries}
\theorembodyfont{\itshape}

\newtheorem{theorem}{Theorem}

\newtheorem{lemma}{Lemma}
\newtheorem{corollary}{Corollary}

\theorembodyfont{\upshape}
{\theoremstyle{break}

}
\theorembodyfont{\upshape}
{\theoremstyle{break}
\newtheorem{remark}{Remark}
}
\newcommand{\qed}{\nopagebreak\par\hspace*{\fill}$\square$\par\vskip2mm}

\newcommand{\id}{\mathop{\mathrm{id}}\nolimits}

\newcommand{\GL}{\mathop{\mathit{GL}}\nolimits}

\newcommand{\tors}{{\mathrm{tors}}}
\newcommand{\Spec}{{\mathrm{Spec}}}
\newcommand{\G}{{\mathcal{G}}}
\newcommand{\Af}{{\mathbb{A}}}

\newcommand{\ind}{{\mathrm{ ind}}}

\title{A Friedlander-Suslin theorem over a noetherian base ring}
\author{Wilberd van der Kallen}
\date{}
\begin{document}
\maketitle
\sloppy
\abstract{
Let $k$ be a noetherian commutative ring and let $G$ be a finite flat group scheme over $k$. Let $G$ act rationally
on a finitely generated commutative $k$-algebra $A$.
We show that the cohomology algebra $H^*(G,A)$ is a finitely generated $k$-algebra. This unifies some earlier results:
If $G$ is a constant group scheme, then it
is a theorem of Evens \cite[Theorem~8.1]{Evens}, and if $k$ is a field of finite characteristic, then it is a theorem of Friedlander and Suslin \cite{Friedlander-Suslin}.
If $k$ is a field of characteristic zero, then there is no higher cohomology, so then it is a theorem of invariant theory.}




\maketitle

\section{Introduction}\label{sec1}

In view of \cite{GoodG} the following theorem will be the key.

\begin{theorem}\label{bounded}
Let $k$ be a noetherian commutative ring and let $G$ be a finite flat group scheme over $k$. 
There is a positive integer $n$  that annihilates $H^i(G,M)$ for all $i>0$ and all $G$-modules $M$.
\end{theorem}

\begin{remark}\label{easyn}
If $k$ does not contain $\mathbb Z$, then one may clearly take $n$ to be the additive order of $1\in k$. 
\end{remark}
\begin{remark}
If $G$ is a constant group scheme, then it is well known that one may take $n$ to be the order of the group
\cite[Theorem 6.5.8]{Weibel}.
(A proof is also implicit in the proof of Theorem \ref{bounded} below.)
\end{remark}

\begin{theorem}[Friedlander-Suslin theorem over noetherian base ring]\label{CFG}
Let $k$ be a noetherian commutative ring and let $G$ be a finite flat group scheme over $k$. 
Let $G$ act rationally
on a finitely generated commutative $k$-algebra $A$.
Then the cohomology algebra $H^*(G,A)$ is a finitely generated $k$-algebra.
\end{theorem}

As usual, cf. \cite[Theorem 1.5]{TvdK}, this implies
\begin{corollary}
If further $M$ is a finitely generated $A$-module with a compatible $G$-action, then $H^*(G,M)$ is a
finitely generated $H^*(G,A)$-module.
\end{corollary}

\begin{remark}
By invariant theory $A^G=H^0(G,A)$ is a finitely generated $k$-algebra \cite[Theorem~8, Theorem~19]{reductive}.
\end{remark}
\begin{remark}
If $G$ is a constant group scheme, then Theorem~\ref{CFG} follows from  \cite[Theorem~8.1]{Evens}. 
The original proof in \cite{Evens} is much more efficient than a proof that relies on \cite{GoodG}.
\end{remark}
\begin{remark}
If $k$ is a field  of finite characteristic, then Theorem~\ref{CFG}  is implicit in \cite{Friedlander-Suslin}. 
(If $G$ is finite reduced, see Evens. If $G$ is finite and connected, take $C = A^G$ in \cite[Theorem~1.5, 1.5.1]{Friedlander-Suslin}. 
If $G$ is not connected, one finishes the argument by following \cite{Evens} as on pages 220--221 of \cite{Friedlander-Suslin}.)

\end{remark}

\begin{remark}
The flatness assumption is essential for this kind of representation theory. One needs it to ensure that taking invariants is
left exact \cite[I~2.10 (4)]{Jantzen} and that the category of comodules is abelian  \cite[I~2.9]{Jantzen}. 
\end{remark}
\section{Discussion}
While the present paper is short,  the full proof of Theorem \ref{CFG} is of book length.
 Presently it uses
\cite{Friedlander-Suslin},  \cite{Touze}, \cite{TvdK}, \cite{FvdK}, \cite{GoodG}, with ideas taken from \cite{Evens}, \cite{BH},
\cite{Srinivas vdK}, \cite{Grosshans book}, \cite{Mathieu G} \dots

\section{Conventions}The coordinate ring $k[G]$ is a Hopf algebra \cite[Part~I, Chapter~2]{Jantzen}. The dual Hopf algebra $M(G)$ is the algebra of measures on $G$. Recall that over a noetherian commutative
ring finite flat modules are finitely generated projective.
Thus both $k[G]$ and $M(G)$ are finitely generated projective $k$-modules.

Following \cite[Part~I, Chapter~8]{Jantzen} we denote by $M(G)^G_\ell$ the $k$-module of left invariant measures. 
Any $G$-module  $V$ may be viewed as a left $M(G)$-module and one has
$M(G)^G_\ell V\subseteq V^G=H^0(G,V)$. But in \cite[Part~I, Chapter~8]{Jantzen} it is assumed that $k$ is a field,  
so we will refer to Pareigis \cite{Pareigis} for some needed facts. 

We say that $G$ acts rationally on a commutative $k$-algebra $A$, if $A$ is also a $G$-module and the multiplication map 
$A\otimes_k A\to A$ is a map of $G$-modules.

An abelian group $L$ is said to have \emph{bounded torsion} if there is a positive integer $n$ so that $nL_\tors=0$, where $L_\tors$
is the torsion subgroup of $L$.
By  \cite[Theorem~10.5]{GoodG} bounded torsion is intimately related to finite generation of cohomology algebras.

\section{Proofs}\label{sec2}

\paragraph{Proof of Theorem~\ref{CFG} assuming Theorem~\ref{bounded}.}
First embed $G$ into some $\GL_N$, as follows. As in \cite[I~2.7]{Jantzen}, let $\rho_r$ denote right regular representation  on $k[G]$. Let $1\in G(k)$ denote the unit element. Observe that $\rho_r(g)(f)(1)=f(g)$
for $g\in G(k)$, $f\in k[G]$. 
It follows that $\rho_r$ is faithful. Further $k[G]$
is a finitely generated projective $k$-module, so there is a $k$-module $Q$ so that $k[G]\oplus Q$ is free, of rank $N$ say.
We let $G$ act trivially on $Q$ and get a faithful action of $G$ on $k[G]\oplus Q\cong k^N$. 

\begin{lemma}
This defines a closed embedding  $G\subset\GL_N$.
\end{lemma}
\paragraph{Proof of Lemma}
Choose algebra generators $f_1,\cdots,f_r$ of $k[G]$. If $R$ is a $k$-algebra then any $g\in G(R)$ is determined by
its values $f_i(g)$. Let $\G(R)$ be the subset of $R^r$ consisting of the $a=(a_1,\cdots,a_r)$ such that $f_i\mapsto a_i$ extends to an algebra homomorphism, denoted $g_a$, from $k[G]$ to $R$.
In other words, $\G$ is the subfunctor associated with the closed embedding $G\subset \Af^r$ given by the $f_i$.
Note that $g_a\in G(R)$ for $a\in \G(R)$.

We seek equations on $\GL_N$ that cut out the image of $G$ in $\GL_N$. It turns out to be more convenient to  cut out an intermediate subscheme $X$.
Let $L\in\GL_N(R)$, viewed
as an $R$-linear automorphism of $R[G]\oplus (R\otimes Q)$. If $L$ is of the form $\rho_r(g)\oplus \id_{R\otimes Q}$ for some $g\in G(R)$, then 
it satisfies
\begin{itemize}\item $L(R[G])=R[G]$, 
\item If $a=(L(f_1)(1),\cdots,L(f_r)(1))$, then $a\in\G(R)$,
\item $g=g_a$.
\end{itemize}
The first two properties define a $k$-closed subscheme $X$ of $\GL_N$ and the last property shows that $G$ is a retract of $X$. Therefore the embedding is closed.

Alternative proof of the Lemma: If $k[G]$ and $Q$ are free $k$-modules, then one may ignore $Q$ and use the proof of \cite[Theorem~3.4]{Waterhouse}, taking $V=A$.
Now use that $G\to \GL_N$ is a closed embedding if there is a cover by affine opens of $\Spec(k)$  over which it is 
a closed embedding.

\qed

We may thus view $\GL_N$ as a group scheme over $k$ with $G$ as $k$-subgroup scheme.
Notice that $\GL_N/G$ is affine  \cite[I~5.5(6)]{Jantzen} 
so that $\ind_G^{\GL_N}$ is exact \cite[Corollary~I~5.13]{Jantzen}. Thus by  \cite[I~4.6]{Jantzen} we may
rewrite $H^*(G,A)$ as $H^*(\GL_N,\ind_G^{\GL_N}(A))$, with $\ind_G^{\GL_N}(A)$ a finitely generated $k$-algebra, by invariant
theory. As $A^G$ is noetherian, it has bounded torsion, and by Theorem~\ref{bounded} $H^{>0}(\GL_N,\ind_G^{\GL_N}(A))=H^{>0}(G,A)$
also has bounded torsion. Theorem~\ref{CFG} now follows from Theorem~10.5 in \cite{GoodG}. 
\qed

Remains to prove Theorem~\ref{bounded}.

\paragraph{Proof of Theorem~\ref{bounded}.}
($\bullet$) Observe that the problem is local in the Zariski
topology on $\Spec(k)$, by the following Lemma.

\begin{lemma}
Let $\cal M$ be a collection of $k$-modules.
Let $f_1,\cdots, f_s \in k$ and let $n_1,\cdots,n_s$ be positive integers, such that $n_i (M\otimes k[1/f_i])_\tors =0$ for 
all $M\in \cal M$ and all $i$.
If the $f_i$ generate the unit ideal, then $n_1\cdots n_s M_\tors =0$ for all $M\in \cal M$.
\end{lemma}
\paragraph{Proof of Lemma} Recall that the $f_i$ generate the unit ideal if and only if the principal open subsets
$D(f_i)= \Spec(k[1/f_i])$ cover $\Spec(k)$.
Take $M\in\cal M$ and $m\in M_\tors$. The annihilator of $n_1\cdots n_s m$ contains a power of $f_i$ for each $i$.
These powers generate the unit ideal.
\qed

Let $H$ be the Hopf algebra $k[G]$. In the notations of Pareigis \cite{Pareigis} we have a rank 1 projective $k$-module $P(H^*)$
that is a direct summand of the $k$-module $H^*=M(G)$ \cite[Lemma~2, Proposition~3]{Pareigis}.  
If that projective module is free, then Pareigis shows that $M(G)^G_\ell$ is a direct summand of $H^*=M(G)$, free of rank one  \cite[Lemma~3]{Pareigis}.
By the observation ($\bullet$) we may and shall assume that  $P(H^*)$ is indeed free. Take a generator  $\psi$ of $M(G)^G_\ell$.
By remark \ref{easyn} we may also assume that $k$ contains $\mathbb Z$, so that tensoring with $\mathbb Q$ does not kill everything.
 
We claim that $\psi(1)$ is now a unit in $k_1=\mathbb Q\otimes k$. It suffices to check this at a geometric point $x=\Spec(F)$
of $\Spec(k_1)$. As $F$ is an algebraically closed field of characteristic zero,
 $G$ is a constant group 
scheme at $x$ by Cartier's Theorem~\cite[11.4, 6.4]{Waterhouse}. The coordinate ring $F[G]$ is now the $F$-algebra of maps from the finite group $G(F)$ to $F$.
Evaluation at an element $g$ of $G(F)$ defines a Dirac measure $\delta_g:F[G]\to F$ and $\psi$ is
a nonzero scalar multiple of the sum $\psi_0=\sum_{g\in G(F)}\delta_g$ of the Dirac measures. Evaluating  $\psi_0$
at $1$ yields the order of $G(F)$, which is indeed invertible in $F$. 

Put $\psi_1=(\psi(1))^{-1}\psi$ in $\mathbb Q\otimes M(G)^G_\ell$. Then $\psi_1(1)=1$ and we conclude that
 $1\in\mathbb Q\otimes k\psi(1)$. 
Then there is $a\in k$  so that $a\psi(1)$ is a positive integer $n$. Put $\phi=a\psi$.
We now observe that $\phi - n$ annihilates $k$ in $k[G]$, and thus annihilates all invariants in $G$-modules.
And for any $G$-module $M$ we have  $\phi(M)\subseteq M^G$ because 
$\phi$ is left invariant.

Consider a short exact sequence of $G$-modules  $$0\to M'\to M\stackrel\pi\to M''\to0.$$ If $m''\in M''^G$, let $m\in M$ be a lift.
One has $0=(\phi-n)m''=\pi\phi(m)-nm''$. As $\phi(m)\in M^G$, we conclude that $n$ annihilates 
  the cokernel 
of $M^G\to M'^G$. Taking $M$ injective, we see that $n$ annihilates $H^1(G,M')$. This applies to arbitrary $G$-modules $M'$.
By dimension shift we get Theorem~\ref{bounded}.
\qed

\begin{remark}
One does not need to use ($\bullet$), because actually $\mathbb Q\otimes M(G)^G_\ell$ maps onto $k_1$,
 even when $M(G)^G_\ell$ is not free over $k$. Indeed,
consider the map $v:\mathbb Q\otimes M(G)^G_\ell\to k_1$ induced by $\chi\mapsto\chi(1):M(G)^G_\ell\to k$. To see
that $v$ is surjective, it suffices again to check at the arbitrary geometric point $x=\Spec(F)$
of $\Spec(k_1)$. 

In fact $\mathbb Q\otimes M(G)^G_\ell$ is always free over $k_1=\mathbb Q\otimes k$.
\end{remark}

\section*{COI}The author declares that he has no conflict of interest.



\end{document}